\documentclass{article}
\usepackage{harvard}
\usepackage{enumerate}

\title{A robust two-level incomplete factorization for (Navier-)Stokes saddle point matrices}

\author{Fred W. Wubs\footnote{Johan Bernoulli Institute of Mathematics and Computing Science, University of
Groningen, P.O.Box 407, 9700 AK Groningen, The Netherlands, Email:
f.w.wubs@rug.nl} ~
and Jonas Thies\footnote{j.thies@rug.nl}
}
\begin{document}
\newcommand{\TODO}[1]{~{\bf TODO:~}#1~}
\newcommand{\Fmat}{$\mathcal{F}$-matrix}
\newcommand{\Fmats}{$\mathcal{F}$-matrices}
\newcommand{\Gmat}{gradient-matrix}
\newcommand{\Gmats}{gradient-matrices}
\newcommand{\qed}{\hfill \setlength{\unitlength}{1mm} \begin{picture}(1,1) \put(0,0){{\framebox(2,2){\strut}}} \end{picture} }
\newtheorem{algorithm}{Algorithm}
\newtheorem{thm}{Theorem}
\newtheorem{lem}{Lemma}
\newtheorem{cor}{Corollary}
\newtheorem{myrule}{Rule}
\maketitle
\begin{abstract}
{We present a new hybrid direct/iterative approach to the solution of a
special class of saddle point matrices arising from the discretization of the steady
incompressible Navier-Stokes equations on an Arakawa C-grid. The two-level method introduced
here has the following properties: (i) it is very robust, even close to the point
where the solution becomes unstable; (ii) a single parameter controls fill and convergence,
making the method
straightforward to use; (iii) the convergence rate is independent of the number of
unknowns; (iv) it can be implemented on distributed memory machines in a natural way;
(v) the matrix on the second level has the same structure and numerical properties
as the original problem, so the method can be applied recursively; (vi) the iteration
takes place in the divergence-free space, so the method qualifies as a `constraint
preconditioner'; (vii) the approach can also be applied to Poisson problems.\par
This work is also relevant for problems in which similar saddle point matrices occur, for
instance when simulating electrical networks, where one has to satisfy Kirchhoff's
conservation law for currents.
}

\end{abstract}

\section{Introduction}

Presently, a typical computational fluid dynamics (CFD) problem may involve millions of
unknowns. They represent velocities and pressures on a grid and are determined by solving a
large sparse linear system of equations.
Robust numerical methods are needed to achieve high fidelity. Therefore one often
resorts to direct (sparse) solvers. In general such a method does not fail as long as the
used precision is enough to handle the posedness of the problem.
However, there are two disadvantages to direct methods. Firstly, the amount of memory
required for the factorization is not linear in the number of unknowns, and when
increasing the problem size one may encounter memory limitations sooner than expected due to
fill generated in the factors. Secondly, all the new elements in the factorization have to
be computed, so that the computing time grows sharply, too. This holds especially for 3D
problems, where the computational complexity of direct methods for partial differential
equations (PDEs) grows with the square of the number of unknowns. \par

For this reason one has to resort to iterative methods for very large applications.
Such methods perform a finite number of iterations to yield an approximate solution.
In theory the accuracy achieved increases with the number of iterations performed.
However, iterative methods are often not robust for complex problems. The iteration process
may stall or diverge and the final approximation may be inaccurate. Furthermore they often
require custom numerics such as preconditioning techniques to be efficient.
\par
The hybrid direct/iterative approach presented here seeks to combine the robustness of
direct solvers with the memory and computational efficiency of iterative methods.
It is based on the direct method recently developed for the Stokes \Fmat~by
\citeasnoun{de_Niet_Wubs_2009}, which has the property that the fill does not increase in
the ``gradient'' and ``divergence'' part of the matrix. To extend this to an incomplete
factorization preconditioner one only has to drop velocity-velocity couplings to limit the
amount of fill. We perform a non-overlapping domain decomposition of the grid, and eliminate
the interior velocities using a direct method. For the remaining variables a
Schur-complement problem has to be solved, which we do by a Krylov subspace method
preconditioned by a novel incomplete factorization preconditioner.
\par
In this paper we start out by giving a survey of previous research in
section~\ref{sec:survey}. In section \ref{sec:directmethod} we will
describe the problem in more detail and review the direct method
developed by \citeasnoun{de_Niet_Wubs_2009}.
In section \ref{sec:ilu} we will introduce the proposed iterative
procedure based on this direct method.
In section \ref{sec:numresults} we present numerical results for a series of increasingly
complex CFD problems: the Poisson, Darcy, Stokes and Navier-Stokes equations.

We conclude in section \ref{sec:Concl} by summarizing the method and results and
giving an outlook on future work.

\section{Survey of previous work\label{sec:survey}}

By \citeasnoun{Benzi_Golub_Liesen_2005} a survey is given of
methods currently in use to solve linear systems from fluid flow
problems. In many cases saddle    point problems can be solved
efficiently by a Krylov subspace iteration \cite{VDVorst_2003}
combined with appropriate preconditioning
\cite{Benzi_Olshanskii,Benzi_Golub_Liesen_2005,de_Niet_Wubs_2006,Elman_2002,Kay_Loghin_Wathen_2002,Elman_2008}.
Often a segregated approach is used, i.e. the velocities are
solved independently from the pressures. This results in inner and
outer iterations, the former for the independent systems, and the
latter to bring the solutions of these systems into balance with
each other. We advocate a fully coupled approach.
\par
The idea of combining direct and iterative methods has been used by \citeasnoun{Henon_2006}
and \citeasnoun{Gaidamour_2008} to solve general sparse linear systems arising from the
discretization of scalar PDEs. As in this paper, they reduce the problem to a
Schur-complement system on the separators of a domain decomposition. The Schur-complement
system is solved iteratively using an ILU factorization. As the structural and numerical
properties are not explicitly preserved, robustness and grid-independence cannot be
ascertained for indefinite problems.

Recently, \citeasnoun{de_Niet_Wubs_2009} proposed a direct method
for the solution of \Fmats, of which the incompressible Stokes
equations on an Arakawa C-grid are a special case.
This special purpose method reduces fill and computation time while preserving the structure
of the equations during the elimination. It still suffers from the weaknesses of direct
methods, but only the number of velocity-velocity couplings increases, not
the number of velocity-pressure couplings.
We believe that a better understanding of the \Fmats~will lead to
generalizations that are of interest to a broader class of indefinite problems and
note that there are applications outside the field of fluid mechanics, e.g. in electronic
circuit simulations \cite{vavasis_1994}, which lead to \Fmats.
\par
For incompressible flow one has to satisfy an incompressibility constraint: the velocity
should be divergence-free. We remark that our iterative technique does not violate the
divergence constraint and therefore belongs to the class of `constraint
preconditioners' \cite{Keller_ea_2000}. For details see section \ref{sec:space}.



\section{${\cal F}$-matrices and the direct solution method \label{sec:directmethod}}

\begin{figure}[t!]
\begin{center}
\setlength{\unitlength}{1.0mm}
\begin{picture}(30,30)(0,0)
\put(0,0){\line(0,1){20}} \put(20,0){\line(0,1){20}}
\put(0,0){\line(1,0){20}} \put(0,20){\line(1,0){20}}
\put(0,10){\multiput(0,0)(10,10){2}{\multiput(0,0)(10,-10){2}{\circle*{1}}}}\put(-3,9){$u$}
\put(21,9){$u$} \put(9,-3){$v$} \put(9,21){$v$}
\put(10,10){\circle*{1}} \put(8,7){$p$}
\end{picture}
\end{center}
\caption{\small Positioning of velocity ($u,v$) and pressure ($p$)
variables in the C-grid.} \label{fig:C-grid}
\end{figure}
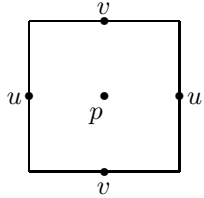

In this paper we study the solution of the equation
\begin{equation}\label{eq:spproblem}Kx = b,\end{equation}
where $K \in R^{(n+m)\times(n+m)}$ $(n\geq m)$ is a saddle point
matrix that has the form
\begin{equation} \label{eq:spmatrix}
K = \left [ \begin{array}{ll} A & B\\ B^T & 0
\end{array}\right ],\end{equation}
with $A\in R^{n\times n}$, $B \in R^{n\times m}$. Special attention is given to a class
of saddle point matrices known as \Fmats. We start out by defining the
gradient matrix in which the \Fmat~ is expressed.
\newtheorem{defi}{Definition}
\begin{defi}\label{dfn:Gmat}
A \Gmat~has at most two nonzero entries per row and its row sum is zero.
\end{defi}
We have chosen the name \Gmat, because this type of matrix
typically results from the discretization of a pressure gradient
in flow equations. It is important to note that the definition
allows a \Gmat~to be non-square. Now we can define the \Fmat.
\begin{defi}\label{dfn:Fmat}
A saddle point matrix (\ref{eq:spmatrix}) is called an \Fmat~ if $A$
is positive definite and $B$ is a \Gmat.
\end{defi}
The definition is due to \citeasnoun{Tuma_2002}.
\Fmats~occur in various fluid flow problems where Arakawa A-grids
(collocated) or C-grids (staggered, see figure~\ref{fig:C-grid})
are used. For example, in \citeasnoun{Arioli_Manzini_2003} the
discretization of Darcy's equation in ground-water flow results in
an \Fmat. They also occur in electronic network simulations \cite{vavasis_1994}.\\

\subsection{The algorithm for the direct approach}\label{sec:algorithm}

Many of the standard algorithms have in common that they compute a
fill-reducing ordering for $K$ and then somehow adapt it to make
it feasible: a factorization is feasible if it does not break down
due to a zero pivot. The delay of elimination (through pivoting) will give an
increase in computing time and may lead to increased fill in the factors.  To preclude this
inefficiency we propose a different approach. Suppose the sets of
all velocities and pressures are denoted by $V$ and $P$,
respectively. The respective elements will be called $V$-nodes
and $P$-nodes. The idea is to first compute an ordering for the
$V$-nodes based on a graph that contains information of the whole
matrix, and then insert the $P$-nodes appropriately. Assume that
we have an elimination ordering on $V$, then we use the following
simple rule to insert $P$-nodes into the ordering:

\begin{myrule} \label{con:uandp} during Gaussian elimination with $K$,
whenever a $V$-node is to be eliminated which is connected to a $P$-node, these nodes
are eliminated together using a $2 \times 2$ pivot.
\end{myrule}

With this rule we get as many $2\times 2$ pivots as
there are $P$-nodes. Only if due to elimination a $V$-node becomes
totally disconnected from $P$ it can be eliminated on its own.

As all $P$-nodes are eliminated together with a $V$-node in
pivots of the form
\[\left(\begin{array}{cc}\alpha & \beta \\\beta & 0\end{array}\right),\]
the factorization is always feasible and additional pivoting is not required.

If we apply this rule to an ordering on $V$ that is constructed as
a fill-reducing ordering for $A$, the resulting ordering for $K$
will not be fill-reducing in general. To ensure that the final
ordering is fill-reducing we have to use information about the
whole matrix, i.e. the fill patterns of $B$ and $B^T$ have
to be taken into account. This is the case if the ordering for $V$ is
fill-reducing for the fill pattern $F(A)\cup F(BB^T)$, where
$F(A)$ denotes the fill pattern of $A$. This graph is an envelope
for the fill that will be created by elimination of the nodes in
$P$. In many cases this will be equal to $F(A+BB^T)$, but to avoid
possible cancellation in the addition we will use the matrix
$F(A)\cup F(BB^T)$.  Summarizing we get the following algorithm:

\begin{algorithm}\label{alg:newalg} To compute a feasible
fill-reducing ordering for the saddle point matrix $K$:
\begin{enumerate}
\setlength{\itemsep}{-0.1cm}
\item Compute a fill-reducing ordering for the $V$-nodes based on $F(A)\cup F(BB^T)$.
\item Insert the $P$-nodes into the ordering according to rule \ref{con:uandp}.
\end{enumerate}
\end{algorithm}
The $P$-nodes (step 2) can be inserted dynamically during Gaussian
elimination, which means that we have to adapt the elimination process.
The elimination is performed using the fill-reducing ordering on $V$ and
applying rule 1. This also takes into account that $V$-nodes initially coupled to
$P$-nodes become decoupled because of cancellation, which is a rather common phenomenon (see
section~\ref{sec:drop}). This is different from just combining
pressures with velocities beforehand (static pivoting).\par

The above method has structure preserving properties which we list
in the theorems below. The first two are taken from
\citeasnoun{de_Niet_Wubs_2009}, where they were proved for symmetric positive definite $A$.
Along the same lines they can be proved for non-symmetric positive definite $A$.

\begin{thm}\label{thmFmat}
If $K$ is an $\mathcal{F}$-matrix, all Schur complements $K^{(l)}$
are $\mathcal{F}$-matrices.
\end{thm}

This means that the $A$ part will remain positive definite and
the $B$ part will have at most 2 entries per row in any step of
the elimination. The latter allows us to keep the $B$ part exact during
the incomplete factorization.

\begin{thm}
The $B$ part in all Schur complements is independent of the size
of the entries in the $A$ part.
\end{thm}

\begin{thm}
 If initially $B$ has entries with magnitude one, then this will remain so during the elimination.
\end{thm}

\begin{thm}
If a $P$-node is not eliminated together with the first $V$-node
it is attached to, the next Schur complement will not be an \Fmat.
\end{thm}

\paragraph{Proof} Consider the matrix in Equation~\ref{eq_gen} in
the next section. It is clear that using only $\alpha$ as pivot
will give a contribution in the zero block. \qed

\vspace*{\baselineskip} Results of the direct method were shown
with AMD \cite{Amestoy_ea_1996} as fill reducing ordering in
\citeasnoun{de_Niet_Wubs_2009}.

\section{Structure preserving incomplete factorization \label{sec:ilu}}

In this section we want to develop an incomplete factorization based on the direct method
described so far. First we will introduce the domain decomposition we use and then we will
illustrate that simply applying a dropping strategy to the $A$ part may not give the
desired result when there are couplings to $P$-nodes. We then proceed to develop a
combination of orthogonal transformations and dropping that leads to grid-independent
convergence, limits fill-in and keeps the divergence constraint intact.

\paragraph{Assumption.}
For this section we will assume that the entries in $B$ have equal magnitude. This is not a restriction because
it can be achieved by scaling the rows of an arbitrary gradient matrix $B$.
If $DB$ gives the desired matrix, our new matrix will be
\[ \left [
 \begin{array}{cc}
  DAD&DB\\
  B^TD&O
 \end{array}
\right ]
\]
Observe that the post-scaling means that the $V$-nodes will be scaled. For Navier-Stokes on a
stretched grid (see section \ref{sec_navstok}) the scaling is such that we get as new unknowns
the fluxes through the control cell boundaries.

\subsection{Domain decomposition}

The first step of the proposed method is to construct a non-overlapping decomposition
of the physical domain into a number of subdomains.
This can be done by applying a graph-partitioning method like Metis
\cite{Karypis_Kumar_1998} or similar libraries to $F(A)\cup F(BB^T)$. Metis has been
tested successfully, but for this paper we use a manual partitioning into equally-sized
square subdomains. (For the Navier-Stokes equations we used a stretched grid, so in that
case they are not square and equally-sized in physical space but in the number of unknowns).
\par

Then we introduce a minimal overlap: two adjacent subdomains share one layer of velocity
nodes, whereas pressure nodes are not shared among
subdomains. Variables belonging to exactly one subdomain are said to be \emph{interior
variables}. Velocities connecting to interior variables in more than one subdomain form
separators of the subdomains they connect to. The separator velocities are complemented by
an arbitrary single $P$-node per subdomain. When eliminating the interior variables in
the next step, this ensures that the subdomain matrix is non-singular (in physical
terms the pressure level inside the subdomain is fixed).
We remark that
\begin{enumerate}[(i)]
\item the domain decomposition can be seen as a Nested Dissection ordering as may be used in
step 1 of Algorithm~\ref{alg:newalg}, stopped at a certain subdomain size (see also
\cite{Toselli_2005} in the paragraph ``Schur Complement Systems'' starting on page 262);

\item we used horizontal and vertical separators as depicted for two domains
in fig.~\ref{fig:ddoms}. A better choice may be to use skew separators
($\pm 45^\circ$), leading to about half the  $V$ nodes on the separator for subdomains of
similar size. Both approaches yield the same number of $V$ nodes with couplings to $P$
nodes in the Schur-complement, and we chose for ease of programming here;

\item we use the decomposition primarily for numerical reasons and the number
of subdomains will typically be much larger
than the number of processors in a parallel computation.

\end{enumerate}

We can now eliminate the interior variables, leading to a Schur-complement problem
for the separator velocities and remaining pressures. The remainder of this section is
devoted to constructing an incomplete factorization preconditioner for this
Schur-complement, so that it can be solved efficiently by a Krylov subspace method.

\subsection{The dropping problem}\label{sec:drop}

Consider the following matrix, which occurs in any elimination step with a $2 \times 2$
pivot:

\begin{equation} \label{eq_gen}
 \left [ \begin{array}{cc|cc}
\alpha & \beta  & a^T & b^T\\
\beta &  0 & \hat b^T &0\\
\hline
a& \hat b& \hat A & \hat B\\
b&  0 &\hat B^T& O
\end{array} \right ].
\end{equation}
When performing the elimination step, a multiple of $\hat b \hat b^T $ is added to $\hat A$.
This does not introduce new fill if $\hat A$ is dense. But if we replaced
$\hat A$ by a sparse matrix by dropping, the matrix would be filled again as $\hat b$ is
typically dense.\par
\begin{figure}[t!]
\begin{center}
\setlength{\unitlength}{1.0mm}
\begin{picture}(50,30)(0,0)
\put(0,0){\line(0,1){20}} \put(20,0){\line(0,1){20}}
\put(40,0){\line(0,1){20}}
\put(0,0){\line(1,0){40}}\put(0,20){\line(1,0){40}}
\put(10,10){\multiput(0,0)(20,0){2}{\circle*{1}}}
\put(-3,9){$u,v$}\put(21,9){$u,v$} \put(41,9){$u,v$}
\put(8,-3){$u,v$}\put(28,-3){$u,v$} \put(8,21){$u,v$}
\put(28,21){$u,v$} \multiput(8,7)(20,0){2}{$p$}
\end{picture}
\end{center}
\caption{\small Velocity separators ($u,v$) and pressure per
domain ($p$) in a 2-domains case.} \label{fig:ddoms}
\end{figure}
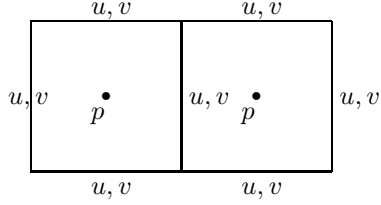

This is a common phenomenon. Consider, for example, the two-domain case in
fig.~\ref{fig:ddoms}. After eliminating the interior variables, many of the $V$-nodes on the
separator are coupled to the two remaining $P$-nodes. Assume that we drop all connections
between the $V$-nodes on the separator, so in the above matrix (\ref{eq_gen}), $\hat A$
is replaced by its diagonal, and $a$ becomes zero; $\hat b$ is a dense vector, $\hat B$ has
an associated dense column with opposite sign, and $b^T$ has a nonzero at the same column
position with sign opposite to that of $\beta$. When eliminating one ``$V$-node $P$-node''
pair, all the $V$-nodes on the separator become detached from $P$ and $\hat A$ becomes
dense.

\par
From the above we learn that we should try
to get more zeros into $\hat b$. Or stated otherwise, we should
try to decouple the $V$-nodes on the separator from the $P$-nodes
as far as possible.

\subsection{Orthogonal operators to decouple $V$- and $P$-nodes} \label{sec:orth}

One idea to get rid of unwanted pressure couplings is to simply drop them. However, the fill
in the $B$-part is already modest and an exact $B$-part is attractive, as discussed
in section~\ref{sec:space}. Fortunately we can do better.
Consider the square domain decomposition (fig. ~\ref{fig:ddoms}), extended periodically so
that every subdomain is bounded by four separators from the neighboring subdomains. The
Schur-complement for the separator velocities and remaining pressures has about the
following form (the $V$-nodes in the corners are neglected here, in practice they form
`separators of the separators' and get a block of their own):
\[ \left [ \begin{array}{cccccc}
A_{11}&B_1&A_{12}& A_{13} & O& O\\
B_1^T & O &B_{21}^T &B_{31}^T&O&O\\
A_{21}&B_{21}&A_{22}&O&A_{24}&B_{22}\\
A_{31}&B_{31}&O&A_{33}&A_{34}&B_{32}\\
O&O&A_{42}&A_{43}&A_{44}&B_{42}\\
O&O&B_{22}^T&B_{32}^T&B_{42}^T&O
\end{array} \right ] \left [
\begin{array}{c}
 v_1\\
 p_1\\
 v_2\\
 v_3\\
 v_4\\
 p_2
\end{array}
\right ] = \left [ \begin{array}{c}
 b_{v_1}\\
 b_{p_1}\\
 b_{v_2}\\
 b_{v_3}\\
 b_{v_4}\\
 b_{p_2}
\end{array}
\right ].
\]
Here $v_1$ contains the $V$-nodes on a certain separator, $p_1$
contains the two $P$-nodes from the adjacent subdomains; $v_2$ and $v_3$ contain the
$V$-nodes from other separators around these subdomains, respectively.
$v_4$ and $p_2$ represent the remaining $V$- and $P$-nodes
in the Schur-complement (separator velocities and pressures not connected to the separator
under consideration).\par

Now $B_1$ only contains two dense columns, equal up to a sign. So
by using an orthogonal  transformation $H$, e.g. a Householder reflection, we can
transform $B_1$ into a matrix with only entries on a certain row, usually
the first. Applying $H$ to the first block row and
column from left and right, respectively, we obtain the following system
(note that the properties of the matrix are preserved by the orthogonal transformation):

\[ \left [ \begin{array}{cccccc}
H^TA_{11}H&H^TB_1&H^TA_{12}&H^TA_{13}&O&O\\
(H^TB_1)^T & O & B_{21}^T &B_{31}^T &O&O\\
A_{21}H&B_{21}&A_{22}&O&A_{24}&B_{22}\\
A_{31}H&B_{31}&O&A_{33}&A_{34}&B_{32}\\
O&O&A_{42}&A_{43}&A_{44}&B_{42}\\
O&O&B_{22}^T&B_{32}^T&B_{42}^T&O
\end{array} \right ] \left [
\begin{array}{c}
 H^Tv_1\\
 p_1\\
 v_2\\
 v_3\\
 v_4\\
 p_2
\end{array}
\right ] = \left [ \begin{array}{c}
 H^Tb_{v_1}\\
 b_{p_1}\\
 b_{v_2}\\
 b_{v_3}\\
 b_{v_4}\\
 b_{p_2}
\end{array}
\right ].
\]
The Householder matrix is a full matrix (though its application is
cheap if its defining form is exploited) and would destroy the
sparsity. However, the matrices $A_{11}$, $A_{12}$ and $A_{13}$ are typically already dense
(see remark below), so not much is lost and we have gained a lot:
we decoupled all but one of the $V$-nodes on the separator from
the $P$-nodes. The decoupled ones can be eliminated on their own now.

\paragraph{Remark 1} The fill of $A_{11}$,
$A_{12}$ and $A_{13}$ depends on the problem at hand. For the 2D Stokes-equations
in the absence of the pressure terms we get two decoupled Poisson
equations for $u$ and $v$. In that case nested dissection gives connections between all the
variables surrounding a domain. So the matrices $A_{11}, A_{12}$, and $A_{13}$ are half full
(no couplings between $u$ and $v$). As most pressures are eliminated with the interior
velocities, the matrices become dense.

\paragraph{Remark 2} In practice, $u$ and $v$ nodes on a separator may connect to the
$P$-nodes with reversed signs. To ensure robustness we apply separate transforms to each
velocity component.

\paragraph{Remark 3} Choosing a Householder transformation may seem arbitrary
and not related to the physics of the problem. We may indeed choose other orthogonal
transformations with the same effect (some alternatives are proposed at the end of section
 \ref{sec:dropstrat}). The key is that one of the columns of $H$ - up to a normalizing
 factor - should
 be the vector $e$ with all entries equal to one. This yields the sum of all the fluxes
 through the interface, so there will be a new variable that represents the entire flux
 through the interface. The other new variables represent fluxes through the interface that
 are on average zero.

\paragraph{Remark 4} Instead of scaling the vectors in $H$ to unit length,
we scale them to the length $m=||e||_2$ of the vector $e$ defining $H$. In that case
the inverse of $H$ is $\frac{1}{m} H^T$.

\paragraph{Remark 5} Although not necessary for the decoupling process, we
also apply an orthogonal transformation to $V$-nodes that are not coupled to a $P$-node in
the first place. This is important for the dropping strategy proposed in the next
section.
\paragraph{}
The situation depicted in eq. \ref{eq_gen} now only occurs once per separator and velocity
component, namely for the $V$-node still coupled to the $P$-nodes.
Because of the transformation $\hat b$ is now zero, and no fill is generated.\par

So far we have not made any approximations, and while we have zeroed out most of the
$V$-node/$P$-node couplings, a dropping strategy has to be applied in the $V$-$V$ part to
get a sparse preconditioner for the Schur-complement. However, the Householder
transformation combined with standard dropping techniques for the SPD case will
generally not lead to grid independent convergence. This requires that the approximation is
spectrally equivalent to the original matrix. We will  consider a new way of dropping in the
next section which has this property.

\subsection{Dropping strategy \label{sec:dropstrat}}

The general idea of the approximation is the following.
We replace the flux through grid cell faces forming a separator
by the combined flux through that separator (see Remark 3 in the previous section). Then we
try to reduce the problem of finding all separator velocities by dropping and elimination to
the related problem of finding the new fluxes (or summed velocities).
This reduced problem can still be understood in terms of conservation of mass and momentum
and its form is very similar to the original problem.\par

Let us consider an orthogonal operator that is more intuitive than the Householder
transformation. Suppose $e$ is a vector with all ones and $C$ is an orthogonal extension of
$e$ such that the length of every column is the same. Define a square matrix \[ H= [ C, e ], \]
which is orthogonal up to a constant factor (see Remark 4 in the previous section).
This operator is applied to the velocity component in normal direction on the separator.
These velocities have the same sign for the connection to the pressure and therefore again
only one row remains in $H^TB_1$. The first component of $H^Tv$ will be the sum of
the components of $v$; we will call this a $V_\Sigma$-node from now on.
To develop some intuition, we first give a simple example of the dropping
strategy which reveals that the resulting reduced problem can be viewed as a coarse
representation of the original one. In section \ref{sec:general} we then perform a more
general analysis.

\subsubsection{Example of dropping}

Consider the familiar tridiagonal matrix with elements $[ -1~2 ~-1]$ on the sub-, main,
and superdiagonal, respectively, which arises when discretizing the 1D Laplace equation
on a regular grid. We premultiply it by the block diagonal matrix $H$
with diagonal blocks
\[ \left [ \begin{array}{cc}
-1&1\\
1&1
\end{array} \right ]
\]
and postmultiply by its transpose (the same matrix here).
Every pair of rows of the transformed matrix has the form
\begin{equation} \label{eq_extrans}
\begin{array}{cccccc}
         1&1&6&0&1&-1 \\
        -1&-1&0&2&1&-1
     \end{array}.
\end{equation}
Next we make an odd-even ordering for the unknowns (equivalent to shifting the
$V_\Sigma$-nodes to the end of the matrix). This new matrix has
the form
\[ \left [ \begin{array}{cc}
A_{11}&A_{12}\\
A_{21}&A_{22}
\end{array} \right ].
\]
The matrix $A_{11}$ is tridiagonal with entries $[1 ~ 6 ~1]$, and $A_{22}$
is tridiagonal with entries $[-1~2~-1]$. So $A_{22}$ is a representation of the original
problem on a twice as coarse grid (up to scaling).
The blocks $A_{12}$ and $A_{21}$ have row sum zero, a typical row of $A_{12}$
being $[1~0~-1]$. We just drop these two blocks and take the remaining part as the
approximation. The fact that both $A_{11}$ and $A_{22}$ are principle submatrices of
the above matrix infers that they both are SPD, so the approximation does not lead to a
singular or indefinite matrix. Note that the elements in the dropped part are quite big and
would not be dropped by traditional drop-by-value strategies.\par

To prove grid-independent convergence when using the resulting matrix as preconditioner, we
have to show that $A_{21} A_{11}^{-1} A_{12} \leq \gamma^2 A_{22}$, for some $\gamma <1$
independent of the size of the matrix
(the case $\gamma=1$ follows directly from  the
positiveness of the Schur-complement of the original problem. For
grid-independence we just need some extra margin).
We can apply Fourier analysis in this constant
coefficient case, which leads to the problem of finding the maximum of
\[  \frac{ \sin(\theta)^2}{(6+2\cos(\theta))\sin(\theta/2)^2} =
  \frac {\cos(\theta/2)^2}{1+ \cos(\theta/2)^2}.
\]
This amounts to finding the maximum of $x/(1+x)$ on [0,1], which is a monotonic
function, so the maximum is 1/2.\par

Another approach is to view the matrix as a sum of ``element'' matrices
$E_i$ and the preconditioner as a sum of $F_i$. Using the Rayleigh quotient,
one can easily show that the condition number of the preconditioned matrix is bounded
if $(x, E_i x)/(x, F_i x)$ is bounded from below and above for $x$ not in the
common null space of $E_i$ and $F_i$ for all $i$ (e.g. see \cite{Axelsson_Larin}).
The singular vector of the transformed matrix (\ref{eq_extrans}) is $[0,1,0,1,0,1,...]^T$,
and for the element matrices of the transformed problem and the approximation we can use
\[
E_i= \left [ \begin{array}{cc|cc}
3&1&1&-1\\
1&1&1&-1\\
\hline
1&1&3&-1\\
-1&-1&-1&1
\end{array}
\right ], ~~
F_i= \left [ \begin{array}{cc|cc}
3&0&1&0\\
0&1&0&-1\\
\hline
1&0&3&0\\
0&-1&0&1
\end{array}
\right ].
\]
Both matrices are nonnegative and the condition number of $E_2^{-1}E_1$ is bounded on
the space orthogonal to $[0,1,0,1]^T$. This approach also reveals that we can
replace $A_{11}$ by any positive diagonal matrix and still have a condition number
independent of the mesh size. This concludes our simple example.

\subsubsection{General analysis \label{sec:general}}

These contemplations suggest that the following lemma and its corollary play a
key role in devising a dropping strategy:
\begin{lem} \label{lemPS} Principal submatrices of an (S)PD-matrix are (S)PD.
\end{lem}
\begin{cor} \label{corPS}
\[ \mbox{If } \left [
\begin{array}{cc}
 A_{11}&A_{12}\\
 A_{21}&A_{22}
\end{array}
\right ]
\mbox{ is (S)PD then }
\left [
\begin{array}{cc}
 A_{11}&O\\
 O&A_{22}
\end{array}
\right ]
\mbox{ is (S)PD.}
\]
\end{cor}
Since we only make approximations in the $A$ part of the matrix $K$, we have the following
lemma.
\begin{lem} \label{lem_bndcond}
If $A$ is SPD, the condition number of the preconditioned $K$ matrix is bounded by the
condition number of the preconditioned $A$, where as preconditioner an SPD approximation of $A$ is used.
\end{lem}
{\bf Proof:} Consider the generalized eigenvalue problem
\begin{equation}
\left [
\begin{array}{cc}
 A-\lambda \tilde{A} & (1-\lambda)B \\
 (1-\lambda) B^T & O
\end{array}
\right ] \left [
\begin{array}{c}
 x_1 \\
 x_2
\end{array}
\right ] = 0,
\end{equation}
where $\tilde A$ denotes an SPD approximation of $A$.
We see that for $\lambda\neq 1$ ($\lambda=1$ is clearly an eigenvalue) we can scale the border by any constant. So the
 eigenvalue problem is in fact an eigenvalue problem restricted to the kernel of the
 divergence (or constraint) operator $B^T$.
 Suppose $Q$ is an orthogonal basis for the kernel of $B^T$, then we have to find the
 eigenvalues of the pencil $(Q^TAQ,Q^T\tilde A Q$. Now
\begin{eqnarray*}
 \lambda_{\mbox{min}} (A,\tilde A)&=&\min_x \frac{(x,Ax)}{(x,\tilde A x)} \leq \min_y \frac{(y,Q^TAQy)}{(y.Q^T\tilde A Qy)} \leq  \frac{(y,Q^TAQy)}{(y.Q^T\tilde A Qy)} \\
&\leq& \max_y \frac{(y,Q^TAQy)}{(y.Q^T\tilde A Qy)} \leq \max_x
\frac{(x,Ax)}{(x,\tilde A x)}=\lambda_{\mbox{max}} (A,\tilde A).
\end{eqnarray*}
 Hence, the eigenvalues of the preconditioned $K$ are bounded by the
 eigenvalues of the preconditioned $A$, which leads to the result. \qed  \\[\baselineskip]
These lemmas set the ground for further reasoning that will lead to grid-independent
convergence. In the remainder of this section we assume that $A$ is symmetric and positive
definite.
%
Let us extend $H$ with an identity for the unknowns that are not transformed and write $H=[H_1,H_2]$, where
\[
 H_1=\left [
\begin{array}{c}
 C\\
0
\end{array}
\right ], ~~
 H_2=\left [
\begin{array}{cc}
e&0\\
0&I
\end{array}
\right ]. ~~
\]
The transformed matrix is given by
\begin{equation} \label{eq_HAH}
H^TAH=
\left [
\begin{array}{cc}
 H^T_1A H_1 &H^T_1 A H_2\\
 H^T_2A H_1& H^T_2 A H_2
\end{array}
\right ].
\end{equation}
Here $ H^T_2 A H_2$ is a Galerkin approximation of $A$ and hence it can be viewed as a
discretization on a coarser grid (in fact it is an aggregation similar to that used by
\citeasnoun{Notay_2010}, albeit Notay applies the aggregation directly to the discretized 
PDE 
whereas we apply it to its Schur complement on the separators). If $A$ is obtained from a
stable discretization of a second-order differential operator, then $H^T_1A H_1$ has a
condition number independent of the mesh size if the dimension of $C$ is fixed (i.e. if the
length of the separator is fixed). We will prove this for a very simple case using finite element theory.
We just consider the operator $\frac{d}{ds} (p(s) \frac{d}{ds} \cdot ) $ with $p >0$ on the 
 interval $(0,d)$ with homogeneous Neumann boundary conditions. Hence the related bilinear 
 form is $a(u,v)=(p u', v')=\int_0^d p u' v'ds$, where we have used the inner product 
 $(u,v)=\int_0^d  u vds$. The norm associated with this inner product is denoted by 
  $||\cdot||$.
Here $u,v$ and $p$ are all functions in the Sobolev space ${\cal H}_1(0,d) $, which 
 consists of all continuous functions that are 
 piecewise differentiable. We will also apply this inner product to vectors of functions,
 which should be read as applying it element by element.
\begin{lem} Let $A=a(V,V)$ and  $M=(V,V)$, where
$V=[\phi_1(s), \phi_2(s), ... \phi_N(s)]$ is a row vector of basis functions in ${\cal H}_1(0,d)$ with the property that
there exists a constant $c$ such that $\frac c{h^2} M - (V',V')$ is nonnegative. If
$c_2 ||u'||^2 \geq a(u,u) \geq c_1 ||u'||^2$, the spectral condition number of
$(H_1^TAH_1)/(H_1^TMH_1)$ is bounded by $(d/h)^2$. 
\end{lem}
{\bf Proof:}
A straightforward substitution of $u=VH_1x$ in the inequality leads to
\[ c_2 ||V'H_1x||^2 \geq a(VH_1x,VH_1x)= (x,H_1^TAH_1x) \geq c_1 ||V'H_1x||^2. \]
Now the minimum of $||f'||/||f||$, where $f$ is an arbitrary ${\cal H}_1$ function
 orthogonal to the constant function is just the eigenfunction of the 1D Laplace operator
 with homogeneous Neumann boundary conditions orthogonal to the constant, which is
 $\cos(\pi s /d)$. So
 \[ \min_x \frac{||V'H_1x||^2}{||VH_1x||^2} \geq \min_f \frac{||f'||^2}{||f||^2} = 
 \left ( \frac \pi d \right )^2. \] 
Hence, the smallest eigenvalue of 
$(H_1^TAH_1)/(H_1^TMH_1)$ is bounded away from zero by $c_1(\frac \pi d )^2$. Now let us try 
to find an upper bound which is less than infinity. This maximum possible is related to the 
highest frequency we can build from $VH_1x$ such that the norm of $V'H_1x$ becomes maximal. 
The shortest wave that can be represented is related to the mesh size $h$. Thus we came to 
 the assumption in the theorem which can be quite easily verified in a special case using 
 for instance the Gershgorin circle theorem. We find that
\[ \frac{||V'H_1x||^2}{||VH_1x||^2} \leq \frac c {h^2} \frac{x^TH_1^TMH_1x}{||VH_1x||^2} = \frac c {h^2},\]
 so the spectral condition number of $H_1^TAH_1/H_1^TMH_1$ is bounded by 
 $\frac{c c_2 d^2}{c_1\pi  h^2}$. \qed  
\par
The constants $c_1$ and $c_2$ are easily determined, we can simply take the minimum and 
maximum of the function $p(s)$ on the interval, respectively.
Although this lemma is based on a second-order differential operator, one could in fact 
find a similar statement for nonnegative operators with pseudo derivative $2 \nu$, where 
$\nu$ may be any positive real number. Such an operator is found, for instance, when 
 writing down the continuous equations at the separators, leading to the so-called 
 Steklov-Poincar\'e operator (see \cite{Toselli_2005}). To return to our discussion, for $d$ 
 in the lemma one could think of the length of the separator. So if $d$ 
 decreases proportional with $h$ when refining the grid, the condition number of  
 $H_1^TAH_1$ is bounded independently of the mesh size assuming we can bound the condition 
 number of $M$ beforehand. 
 The latter matrix is usually strictly diagonally dominant. So simply applying Gershgorin's 
 theorem makes the assumption valid.

Now assume we have the following strengthened Cauchy-Schwarz inequality \cite[section 9.1]{Axelsson94}
\begin{equation} \label{eq_CBS} |x^T H^T_1 A H_2 y | \leq \gamma \{ (x^T H^T_1 A H_1 x) (y^T H^T_2  A H_2 y ) \}^{\frac12}
\end{equation}
holding independently of the mesh size. In our case $H_1^TH_2=0$, and if the
columns of $H_1$ or $H_2$ span an invariant subspace of $A$, then also  $H^T_1 A H_2=0$,
 hence $\gamma=0$.  The latter is only approximately the case here, so we will find some
 $\gamma <1$. Lemma 9.2 from \cite{Axelsson94} states that
\begin{equation}
\label{eq_indep}
H^T_2A H_1 (H^T_1A H_1)^{-1}  H^T_1 A H_2 \leq \gamma^2  H^T_2 A H_2,
\end{equation}
where the inequality should be understood in the sense that the sum of the left-hand side 
and
the right-hand side gives a non-negative matrix.
For ease of notation we write the transformed matrix (\ref{eq_HAH}) as
\[
\left [
\begin{array}{cc}
 A_{11}&A_{12}\\
 A_{21}&A_{22}
\end{array}
\right ].
\]
In this notation the above property reads
$A_{21} A_{11}^{-1} A_{12} \leq \gamma^2 A_{22}$. Now the preconditioner obtained
by dropping $A_{21}$ and $A_{12}$ is SPD according to Corollary~\ref{corPS}. The
eigenvalues of the preconditioned $A$ matrix can be found from the following generalized
 eigenvalue problem.
\[
\label{eqn:geneigprob}
\left [
\begin{array}{cc}
 (1-\lambda)A_{11}&-\lambda A_{12}\\
 -\lambda A_{21}&(1-\lambda)A_{22}
\end{array}
\right ]
 \left [
\begin{array}{c}
 x_{1}\\
 x_{2}
\end{array}
\right ] = 0,
\]
which leads to $((1-\lambda)^2 A_{22}-\lambda^2 A_{21} A_{11}^{-1} A_{12})x_2=0$ for
 $\lambda \neq 1$. Combined with the previous, we can only find eigenvalues for
  $(1-\lambda)^2 < (\lambda \gamma)^2$, so $1 - \lambda \gamma < \lambda < 1 + \lambda
   \gamma  $. 
   So we find that the  condition number of the preconditioned matrix is less than 
   $(1+\gamma)/(1-\gamma)$, where $\gamma$ is
 independent of the mesh size. Using Lemma~\ref{lem_bndcond} we find the main result of this paper.
\begin{thm}\label{mainresult}
If a strengthened Schwarz-inequality (\ref{eq_CBS}) holds for $0 \leq \gamma <1$ independent
of the mesh size, then we have convergence independent of the mesh-size when the dropping
process as discussed above is applied. The condition number of the preconditioned $K$ matrix is
 bounded by $(1+\gamma)/(1-\gamma)$.
\end{thm}


The situation above remains the same if we apply the transformation to all separators at
 once. After the transformation, only the unknowns associated with $A_{22}$ are coupled to
 pressures. We may still have couplings between various separators in $A_{11}$, but the
 condition number of that matrix is independent of the mesh size. To lower the computational
 cost we also drop couplings between separators in $A_{11}$. We conclude this section by a
 number of remarks concerning the dropping strategy.

\paragraph{Scalar equations.}
The reader may have noticed that in this section we hardly mentioned the pressure.
In fact, the combination of orthogonal transformations and dropping may also be
applied to the pure diffusion problem. In section \ref{sec:numresults}
we will start out by showing numerical results for the scalar Poisson equation.

\paragraph{The nonsymmetric case.}
One may ask how much of the above can be generalized to the
nonsymmetric case (for instance the Navier-Stokes equations). Assume
the nonsymmetric matrix $A$ is positive definite (PD), i.e. $(x,Ax)>0$ for any non-trivial
$x$. Then the Schur complement is PD and the orthogonal transformation does not destroy that
property. Since all principle submatrices of a PD matrix are PD, the approximation will be
PD. So the factorization will not break down. To say something about the condition number of
the preconditioned matrix is more difficult. For a mild deviation from symmetry we
expect the same behavior as for the symmetric case. However, the numerical results for the
Navier-Stokes equations at relatively high Reynolds-numbers indicate that the method works
very well even for highly non-symmetric matrices.

\paragraph{Numerical stability.}
In traditional lumping, only possible for M-matrices, one simply lumps a
 coefficient on the diagonal. This means that a nonnegative matrix is subtracted.
 Eijkhout \cite{Victor} showed already in the nineties that this may give a zero on the
 diagonal. This is easy to preclude by simply not allowing the diagonal to become zero.
 What is much harder to prevent is the occurrence of independent systems in the
 preconditioner, some of which may be singular. This easily
 occurs in anisotropic problems. The proposed dropping does not suffer from these
 problems.

\paragraph{Alternatives for $\mathbf{H}$.}
Finally we propose a simple orthogonal extension to $e$ in order to form $H$. Let $m$ be
 the order of $H$ and note that $[1,-1,0,\cdots,0]^T,$ $[1,1,-2 , 0,\cdots,0]^T, \cdots,$
 $[1,\cdots,1,-(m-1)]^T,e$ are all orthogonal. They can be used for the extension after
 a proper scaling to the length of $e$. The application of this operator can be implemented
 by keeping a partial sum. In this way about $2m$ additions of rows of the matrix it is
 applied to are needed. The Householder transform has a similar operation count.
 One may ask whether alternative choices for $C$ in $H_1$ influence the convergence. This
 is not the case. We can replace $H_1$ by $H_1Q$. For arbitrary orthogonal matrices $Q$
 this has no influence on (\ref{eq_CBS},\ref{eq_indep}) and the
 following analysis.


\subsection{Iteration in the kernel of $B^T$}\label{sec:space}

Since the fill of the $B$ part remains at most 2 per row during
the whole process, we will not drop there. This means that the $B$
matrix is exact in the factorization, and with appropriate dropping
(such as the strategy introduced in the previous section),
the eigenvalues of the preconditioned matrix will all be positive
and real. Still, we cannot directly apply the preconditioned
conjugate gradient method since for that both original and
preconditioner must be positive definite in the Krylov subspace.
We can enforce this condition by building the Krylov subspace $\mathcal{K}(\tilde K^{-1}K,x)$ on
a starting solution $x$ that satisfies the constraint. In exact arithmetic
$\mathcal{K}$ then remains in the kernel of $B^T$. In practice, accumulation of round-off
errors will undermine this property.\par
This problem is often encountered in the field of constraint optimization, and
\citeasnoun{Gould_ea_2001} have developed a variant of the conjugate gradient method,
Projected Preconditioned CG (PPCG), which can be used for the Stokes problem.
There are various ways to find a particular solution of $B^T v=b_2$, one of which is solving
the system once, replacing $K$ by the preconditioner.\par

For the Navier-Stokes equations one could devise a Projected Preconditioned FOM method, as
long as the eigenvalues of the preconditioned matrix are in the right half plane, but
for the results shown in section~\ref{sec_navstok} we simply used MATLAB's gmres.

\subsection{Program structure}\label{sec:alg}

Before looking at numerical results, let us review the complete algorithm and
remark on some implementation issues. The main structure of the program is as
follows
\begin{enumerate}
 \item Perform a domain decomposition on $F(A) \cup F(B)F(B)^T$. We just make a rectangular decomposition of the domain here.
\item Group the variables into subdomain variables and separator variables
(velocities connecting to variables in more than one subdomain).
All pressures are treated as subdomain variables at this stage.
\item \label{groupingstep} Group the separator variables according to variable type (i.e. $u$, $v$) and the subdomains they have connections to.
Thus, we will get a group of $u$-velocities connecting to variables on subdomains 1 and 2, for instance.
\item In the corners of subdomains a complete conservation cell (see fig.~\ref{fig:C-grid})
can occur on the separators.
 This would lead to a singularity in step \ref{elimsdvars}. The velocities making
 up such a cell are flagged '$V_\Sigma$'-nodes (cf. step \ref{pickVsums}). Both these
 $V_\Sigma$-nodes and the $P$-node in the cell will be retained in the Schur-complement.
\item Pick for every domain a $P$-node to be kept in the reduction, and shift these to the
end of the ordering (i.e. retain them in the Schur-complement).
\item \label{elimsdvars} Eliminate all interior variables of the subdomains and construct
the Schur complement system for the velocity separators and the selected pressures.
\item Perform the transformation on each separator group identified in step
\ref{groupingstep}.
\item \label{pickVsums} Identify $V_\Sigma$ nodes (separator velocities that still connect
to two pressures) and put them at the end of the ordering, just before the remaining
pressure nodes.
\item Drop  all connections between non-$V_\Sigma$ nodes and $V_\Sigma$ nodes, and between
non-$V_\Sigma$  nodes in different separator groups. The resulting matrix
is block-diagonal with the `reduced Schur-complement' in the lower right corner.
\item Iterate on the Schur complement using the matrix of the previous step as
preconditioner. This preconditioner is easily applied using LU decompositions
of all non-$V_\Sigma$ blocks and the reduced system.
\end{enumerate}

In three space dimensions, step \ref{groupingstep} is implemented by first numbering the faces, then the edges and then the corners of the box-shaped
subdomains. We note that this is a special case of the hierarchical interface decomposition
(HID) used by \citeasnoun{Henon_2006} and \citeasnoun{Gaidamour_2008}.

\subsection{Computational complexity}

We will now discuss the complexity of the algorithm, implemented as discussed in the
previous section. We assume that a direct method with optimal complexity is used for the
solution of the relevant linear systems, so in 3D if the number of unknowns is
$\mathcal{O}(N)$, the work is $\mathcal{O}(N^2)$, as with Nested Dissection.
For the 3D (Navier-)Stokes equations, we have $N=\mathcal{O}(n^3)$ unknowns, where
$n$ is the number of grid cells in one space dimension. We keep the subdomain size
constant and denote the number of unknowns per subdomain by $S=\mathcal{O}(s^3)$
(here $s$ is the fixed separator length). Hence, there will be $N/S$ subdomains.
Per domain there will be $\mathcal{O}(s^2)$
non-$V_\Sigma$- and $\mathcal{O}(1)$ $V_\Sigma$-nodes.
Per domain the amount of work required is as follows:
\begin{enumerate}
\item $\mathcal{O}(S^2)$ for the subdomain elimination;
\item transformation on faces with $H$: $\mathcal{O}(s^4)$;
\item factorization of non-$V_\Sigma$ nodes: $\mathcal{O}((s^2)^3)=O(S^2)$.
\end{enumerate}
The total over all domains is $\mathcal{O}(N/S) \mathcal{O} (S^2)= \mathcal{O} (NS)$,
 so in this part the number of operations decreases linearly with $S$ (e.g. by a factor 8
 if $s$ is halved).\par
 The solution of the reduced problem ($V_\Sigma$-nodes) requires
 $\mathcal{O}((N/S)^2)$ operations. Here doubling $s$ will decrease the work by
 a factor 64. So in total the work per
 iteration is $\mathcal{O} (NS)+\mathcal{O}( (N/S)^2)$. The number of iterations is constant
  for $S$ constant. There is, however, a positive dependence on $S$ as we may expect. In the
  next section we will observe that the number of iterations is proportional to $\log(S)$. So
  if we double $s$, a fixed amount of iterations is added.  \par

It is clear that if we solved
the reduced problem iteratively by applying our method recursively until the problem
has a fixed grid-independent size, the overall complexity would be
$\mathrm{log}(S)\mathcal{O}(NS)$.

\section{Numerical experiments \label{sec:numresults}}

In this section we will demonstrate the performance of the new solver by applying it
to a series of increasingly complex problems relevant to computational fluid dynamics. For
each problem we first keep the subdomain size constant while refining the mesh.
As discussed in the previous section, the complexity of the algorithm will then be linear in
the number of unknowns except when solving the reduced Schur complement: the operations
required to factor a single subdomain matrix stays the same
and the number of subdomains increases linearly with the grid size. Furthermore, both size
and connectivity pattern of the separators remain the same so the amount of work per
separator remains constant while the number of separators increases linearly, too.\par
The second experiment will be to fix the grid size and vary the subdomain size (i.e. the
number of subdomains). The expectation here is that due to fill-in the bulk of the work load
shifts from the Schur-complement towards the subdomain factorization as the size of the
subdomains is increased.
\par

For each experiment, the following data is displayed:
\begin{itemize}
\item $n_x$ - the grid size is $n_x \times n_x$ ($n_x \times n_x \times n_x$) in 2D (3D),
respectively.
\item $s_x$ - the subdomain size is $s_x \times s_x$ ($s_x \times s_x \times s_x$) in 2D
(3D), respectively.
\item $\mathrm{N}$ - number of unknowns (size of the saddle point matrix),
\item nnz - number of nonzeros in original matrix,
\item $\mathrm{N_S}$ - number of unknowns on the separators and remaining p's (size of the
Schur-complement),
\item $\mathrm{n}$ - number of $V_\Sigma's$ and remaining p's (size of reduced
Schur-complement),
\item iter - number of CG iterations performed on the Schur-complement to reduce the
residual norm by $\mathrm{1/tol}=10^8$,
\item fill 1 - grid-independent part of relative fill-in (number of nonzeros in the
solver divided by number of nonzeros in original matrix). The grid-independent portion
consists of
\begin{itemize}
\item a) fill-in generated while factoring the subdomain matrices
\item b) fill-in generated while constructing the Schur-complement
\item c) fill-in generated while factoring the separator-blocks of the preconditioner
\end{itemize}
\item fill 2 - grid-dependent part of relative fill-in, generated when factoring the
$n\times n$-dimensional reduced Schur-complement.
\item $\kappa$ - condition estimate of the preconditioned Schur-complement: fraction
of the largest and smallest eigenvalue (by magnitude) of the generalized eigenvalue problem
$\mathrm{S}x+\lambda \mathrm{M}x = 0$, where $\mathrm{S}$ is the Schur-complement, and
$\mathrm{M}$ the preconditioner used. We use approximations to the actual eigenvalues
computed by MATLAB's `eigs' command (Not all tables contain this value).
\end{itemize}
\paragraph{Remark} the fill listed under fill 1 b) can be avoided by not explicitly constructing
the Schur-complement. The fill listed as 'fill 2' grows with increasing grid size, but it
can be made grid-independent by solving $S_2$ iteratively, too (i.e. by applying our method
recursively).\\[\baselineskip]

We do not show plots of the convergence behavior. Since all the results are obtained by CG the
convergence is, apart from the first few digits gained, completely regular, which shows that
the eigenvalues, except for a few outliers at the beginning, appear in a cluster.
The relatively stringent convergence tolerance of 8 digits ensures that the overall
convergence behavior does not strongly depend on the choice of the initial
vector. Choosing a smaller tolerance results in stagnation for some of the
tests below because the conditioning of the matrix doesn't allow for more accurate
solutions.

The general behavior we observe in the second experiment
 is that the number of iterations scales with $\log(s_x)$, where $s_x$ is the separator
 length. So doubling the separator length means an increase of the number of iterations by a
 constant amount.

\subsection{The Poisson equation}

We first investigate Poisson's equation, discretized using second order central differences
on a regular structured grid (standard 5-point and 7-point stencils in 2D and 3D,
respectively). This is an important case as solving Poisson's equation is
central to most CFD problems, for instance to determine the pressure in explicit time
stepping algorithms. Tables~\ref{ta_Pois2D_1} and \ref{ta_Pois2D_2}  show the 2D results.
The first shows the dependence on grid refinement and the latter the influence of the domain sizes.
Similar results for the 3D case are shown in tables~\ref{ta_Pois3D_1} and \ref{ta_Pois3D_2}.

\begin{table}[!ht]
 \begin{tabular}{|c|r|r|r|r|r|r|r|r|}
\hline
$n_x$ & N & nnz & $\mathrm{N_S}$ & n & iter & fill 1 & fill 2 & $\kappa$ \\
\hline
32 & 1 024 & 5 112 &  240 & 48 & 21 & 5.53 & 0.20 & 7.04 \\
64 & 4 096 & 20 472 &  960 &  192 & 21 & 5.52 & 0.39 & 7.04 \\
 128 & 16 384 & 81 912 & 3 840 &  768 & 21 & 5.52 & 0.68 & 7.04 \\
 256 & 65 536 &  327 672 & 15 360 & 3 072 & 21 & 5.52 & 1.03 & 7.04 \\
 512 &  262 144 & 1 310 712 & 61 440 & 12 288 & 21 & 5.52 & 1.59 & 7.04 \\
1 024 & 1 048 576 & 5 242 872 &  245 760 & 49 152 & 21 & 5.52 & 2.20 & 7.04 \\
\hline
\end{tabular}
\caption{2D Poisson-equation - grid refinement, subdomain size $s_x=8$.
\label{ta_Pois2D_1}}
\end{table}

\begin{table}[!ht]
 \begin{tabular}{|c|r|r|r|r|r|r|r|r|}
\hline
$s_x$ & N & nnz & $\mathrm{N_S}$ & n & iter & fill 1 & fill 2 & $\kappa$ \\
\hline
4  & 1 048 576 & 5 242 872 &  458 752 &  196 608 & 16 & 2.01 & 11.5 & 4.00 \\
8  & 1 048 576 & 5 242 872 &  245 760 & 49 152 & 21 & 5.52 & 2.29 & 7.04 \\
16 & 1 048 576 & 5 242 872 &  126 976 & 12 288 & 27 & 9.84 & 0.39 & 11.2 \\
32 & 1 048 576 & 5 242 872 & 64 512 & 3 072 & 32 & 13.8 & 0.063 & 16.5 \\
\hline
\end{tabular} \\

\caption{2D Poisson-equation - increasing subdomain size, grid-size $n_x=1 024$\label{ta_Pois2D_2}}
\end{table}

\begin{table}[!ht]
 \begin{tabular}{|c|r|r|r|r|r|r|r|r|r|}
\hline
$n_x$ & N & nnz & $\mathrm{N_S}$ & n & iter & fill 1 & fill 2 & $\kappa$ \\
\hline
16 & 4 096 & 28 660 & 1 352 & 56 & 24 & 29.7 & 0.064 & 10.1 \\
32 & 32 768 &  229 364 & 10 816 &  448 & 25 & 29.0 & 0.36 & 10.2 \\
64 &  262 144 & 1 834 996 & 86 528 & 3 584 & 25 & 29.0 & 1.53 & -  \\
\hline
\end{tabular}
\caption{3D Poisson-equation - grid refinement, subdomain size $s_x=8$ \label{ta_Pois3D_1}}
\end{table}

\begin{table}[!ht]
\begin{tabular}{|c|r|r|r|r|r|r|r|r|}
\hline
$s_x$ & N & nnz & $\mathrm{N_S}$ & n & iter & fill 1 & fill 2 & $\kappa^\star$\\
\hline
4 &  262 144 & 1 834 996 &  151 552 & 28 672 & 19 & 3.68 & 52.0 & 5.75 \\
8 &  262 144 & 1 834 996 & 86 528 & 3 584 & 25 & 29.0 & 1.5 & 10.2 \\
16 &  262 144 & 1 834 996 & 46 144 &  448 & 30 & 116.2 & 0.045 & 16.7 \\
\hline
\end{tabular}
\caption{3D Poisson-equation - increasing subdomain size, grid size $n_x=64$.\newline
$\star$ Computed at $n_x=32$.
\label{ta_Pois3D_2}}
\end{table}

\subsection{Darcy's law}

For flows in porous media one often has to solve the Darcy problem, where $A$ is just a
diagonal matrix. One approach is to eliminate the velocities, which leads  to a Poisson
equation. Care has to be taken when calculating the velocities, because
the gradient operator has to be applied to the pressure. In this numerical differentiation
of the pressure field, round-off errors may be amplified too much to obtain an accurate
solution. Therefore, Darcy's problem is often solved in primitive form.
Tables \ref{ta_Darcy2D_1} through \ref{ta_Darcy3D_2} show the numerical results for Darcy's
law in two and three space dimensions.
\begin{table}[!ht]
 \begin{tabular}{|c|r|r|r|r|r|r|r|r|}
\hline
$n_x$ & N & nnz & $\mathrm{N_S}$ & n & iter & fill 1 & fill 2 & $\kappa$ \\
\hline
16 &  736 & 2 400 & 65 & 17 & 16 & 5.53 & 0.061 & 3.77 \\
32 & 3 008 & 9 920 &  385 &  109 & 25 & 6.29 & 0.24 & 10.8\\
64 & 12 160 & 40 320 & 1 793 &  533 & 26 & 6.65 & 0.49 & 12.2\\
 128 & 48 896 &  162 560 & 7 681 & 2 341 & 26 & 6.82 & 1.00 & 12.6\\
 256 &  196 096 &  652 800 & 31 745 & 9 797 & 26 & 6.91 & 1.69 & 12.6\\
 512 &  785 408 & 2 616 320 &  129 025 & 40 069 & 26 & 6.95 & 2.64 & 12.7 \\
1 024 & 3 143 680 & 10 475 520 &  520 193 &  162 053 & 26 & 6.97 & 3.58 & - \\
\hline
\end{tabular}
\caption{2D Darcy-equation - grid refinement, subdomain size $s_x=8$\label{ta_Darcy2D_1}}
\end{table}

\begin{table}[!ht]
 \begin{tabular}{|c|r|r|r|r|r|r|r|r|}
\hline
$s_x$ & N & nnz & $\mathrm{N_S}$ & n & iter & fill 1 & fill 2 &$\kappa$ \\
\hline
8 & 3 143 680 & 10 475 520 &  520 193 &  162 053 & 26 & 6.97 & 3.58 & 12.7\\
16 & 3 143 680 & 10 475 520 &  258 049 & 40 069 & 29 & 11.1 & 0.66 & 17.6 \\
\hline
\end{tabular}
\caption{2D Darcy-equation - increasing subdomain size, grid size $n_x=512$\label{ta_Darcy2D_2}}
\end{table}

\begin{table}[!ht]
 \begin{tabular}{|c|r|r|r|r|r|r|r|r|r|}
\hline
$n_x$ & N & nnz & $\mathrm{N_S}$ & n & iter & fill 1 & fill 2 & $\kappa$ \\
\hline
8 & 1 856 & 6 720 &  492 &  171 & 34 & 10.8 & 1.28 & 14.0 \\
16 & 15 616 & 57 600 & 5 878 & 2 683 & 36 & 10.2 & 17.6 & 15.3 \\
32 &  128 000 &  476 160 & 54 762 & 27 819 & 36 & 9.73 & 87.7 & 15.4 \\
40 &  251 200 &  936 000 &  109 972 & 56 971 & 36 & 9.65 & 167. & - \\
\hline
\end{tabular}

\caption{3D Darcy-equation - grid refinement, subdomain size $s_x=4$ \label{ta_Darcy3D_1}}
\end{table}

\begin{table}[!ht]

 \begin{tabular}{|c|r|r|r|r|r|r|r|r|r|}
\hline
$s_x$ & N & nnz & $\mathrm{N_S}$ & n & iter & fill 1 & fill 2 & $\kappa^\star$ \\
\hline
4 &  251 200 &  936 000 &  109 972 & 56 971 & 36 & 9.65 & 167. & 15.4 \\
8 &  251 200 &  936 000 & 53 037 & 11 601 & 39 & 50.2 & 16.7 & 18.3 \\
\hline
\end{tabular}
\caption{3D Darcy-equation - increasing subdomain size, grid size 
$n_x=40$. \newline
$\star$ Computed at $n_x=32$.}
\label{ta_Darcy3D_2}
\end{table}

\subsection{A Stokes problem}

\label{sec:discr} The problem is a two-dimensional Stokes equation
on the unit square
\begin{equation}\left. \begin{array}{rcl} -\nu \Delta \mathbf{u} + \nabla p & = & 0 ~, \\
\nabla \cdot \mathbf{u} & = & 0 ~, \end{array} \right\}\end{equation} where
$\mathbf{u}(x,y)$ is the velocity field and $p(x,y)$ the pressure field;
the parameter $\nu$ controls the amount of viscosity. We can get rid of the parameter $\nu$
by defining a new pressure variable $\bar p = p/\nu$. If the first equation is divided by
$\nu$, we can substitute $p$ by $\bar p$ and the parameter $\nu$ is gone. So we may assume
that $\nu = 1$.

These equations are discretized on a uniform staggered grid (a
C-grid, see fig.~\ref{fig:C-grid}) which results in an \Fmat. It
is singular because the pressure field is determined up to a
constant.

For the Stokes problem the matrix $B^T$ represents the
discrete divergence operator. Consequently, we call the kernel of
this matrix the divergence free space.
As a solution of this problem we choose a random vector
in the divergence free space. So the right-hand side of the
divergence equation is zero in our case.

We start off the iteration with the zero vector (which is
trivially in the divergence free space) and therefore we can use
the projected conjugate gradient method (see section~\ref{sec:space}).
Results are summarized in tables \ref{ta_Stokes2D_1} through \ref{ta_Stokes3D_2}.

\begin{table}[!ht]
 \begin{tabular}{|c|r|r|r|r|r|r|r|r|r|}
\hline
$n_x$ & N & nnz & $\mathrm{N_S}$ & n & iter & fill 1 & fill 2 & $\kappa$ \\
\hline
16 &  736 & 4 196 & 65 & 17 & 18 & 7.79 & 0.057 & 4.93 \\
32 & 3 008 & 17 604 &  385 &  109 & 27 & 8.39 & 0.25 & 12.8 \\
64 & 12 160 & 72 068 & 1 793 &  533 & 31 & 8.68 & 0.65 & 13.8 \\
 128 & 48 896 &  291 588 & 7 681 & 2 341 & 31 & 8.72 & 1.33 & 14.2 \\
 256 &  196 096 & 1 172 996 & 31 745 & 9 797 & 31 & 8.70 & 2.40 & 14.6 \\
 512 &  785 408 & 4 705 284 &  129 025 & 40 069 & 31 & 8.60 & 3.83 & 15.0 \\
\hline
\end{tabular} \\

\caption{2D Stokes-equation - grid refinement, subdomain size $s_x=8$\label{ta_Stokes2D_1}}
\end{table}

\begin{table}[!ht]
 \begin{tabular}{|c|r|r|r|r|r|r|r|r|r|}
\hline
$s_x$ & N & nnz & $\mathrm{N_S}$ & n & iter & fill 1 & fill 2 & $\kappa$ \\
\hline
 4 &  785 408 & 4 705 284 &  260 097 &  162 053 & 24 & 3.65 & 20.0 & 9.6 \\
 8 &  785 408 & 4 705 284 &  129 025 & 40 069 & 31 & 8.60 & 3.83 & 15.0 \\
16 &  785 408 & 4 705 284 & 63 489 & 9 797 & 38 & 15.7 & 0.60 & 21.9 \\
\hline
\end{tabular} \\

\caption{2D Stokes-equation - increasing subdomain size, grid size $n_x=512$\label{ta_Stokes2D_2}}
\end{table}

\begin{table}[!ht]
 \begin{tabular}{|c|r|r|r|r|r|r|r|r|r|}
\hline
$n_x$ & N & nnz & $\mathrm{N_S}$ & n & iter & fill 1 & fill 2 & $\kappa$ \\
\hline
8 & 1 856 & 13 728 &  492 &  171 & 34 & 13.9 & 1.20 & 16.6 \\
16 & 15 616 &  122 304 & 5 878 & 2 683 & 41 & 12.5 & 16.4 & 23.8 \\
32 &  128 000 & 1 029 504 & 54 762 & 27 819 & 43 & 11.5 & 103. & 27.1 \\
40 &  251 200 & 2 030 880 &  109 972 & 56 971 & 43 & 11.3 & 168. & - \\
\hline
\end{tabular}

\caption{3D Stokes-equation - grid refinement, subdomain size $s_x=4$ \label{ta_Stokes3D_1}}
\end{table}

\begin{table}[!ht]
 \begin{tabular}{|c|r|r|r|r|r|r|r|r|r|}
\hline
$s_x$ & N & nnz & $\mathrm{N_S}$ & n & iter & fill 1 & fill 2 & $\kappa^\star$\\
\hline
 4 &  251 200 & 2 030 880 &  109 972 & 56 971 & 43 & 11.3 & 167. & 27.1 \\
 8 &  251 200 & 2 030 880 & 53 037 & 11 601 & 49 & 65.8 & 12.1 & 39.1 \\
\hline
\end{tabular} \\

\caption{3D Stokes-equation - increasing subdomain size, grid size 
$n_x=40$.\newline
$\star$ Computed at $n_x=32$.
\label{ta_Stokes3D_2}}
\end{table}

\subsection{Incompressible flow in a lid-driven cavity \label{sec_navstok}}

As test problem for the Navier-Stokes equations we use the lid driven cavity. In
\cite{Tiesinga_2002} this problem was studied near the transition point from steady to
 transient flow. The stability of steady and periodic solutions
was investigated using the Newton-Picard method \cite{Lust_1999} with
  the $\theta$-method for time stepping (with $\theta$ slightly larger than $0.5$ in order
  to damp high-frequency modes which would otherwise show up as spurious eigenvalues
  near the imaginary axis). The linear systems that have to be solved have a slightly
  increased diagonal, which improves the conditioning somewhat.
  The MRILU preconditioner \cite{BottaWubs} used at the time
  converged slowly and not at a grid-independent rate. 
 In a recent review \cite{Elman_2008}, the performance of a number of block multi-level 
 preconditioners is investigated for the steady problem for Reynolds numbers up to 1000. 
 These methods also solve the coupled equations, but perform inner iterations on the 
 velocity and pressure part separately and hence require many parameters to be tuned. 
 Below we demonstrate robust, grid-independent convergence for the driven cavity problem
 at Reynolds-numbers of up to 8000.\par
 
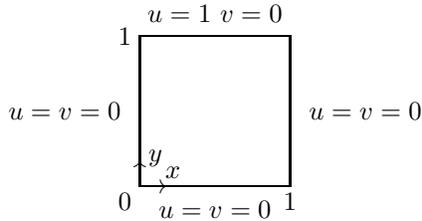
\begin{figure}[bt]
\setlength{\unitlength}{2cm}
\begin{center}
\begin{picture}(1,1.2)(0,0)
\put(0,0){\line(0,1){1}}
\put(1,0){\line(0,1){1}}
\put(0,1){\line(1,0){1}}
\put(0,0){\line(1,0){1}}
\put(0.5,1.15){\makebox(0,0)[c]{$u=1~v=0$}}
\put(0.5,-0.15){\makebox(0,0)[c]{$u=v=0$}}
\put(-0.5,0.5){\makebox(0,0)[c]{$u=v=0$}}
\put(1.5,0.5){\makebox(0,0)[c]{$u=v=0$}}
\put(0.1,-0.01){\makebox(0,0)[c]{$\rightarrow$}}
\put(0,0.1){\makebox(0,0)[c]{$\uparrow$}}
\put(0.22,0.09){\makebox(0,0)[c]{$x$}}
\put(0.105,0.19){\makebox(0,0)[c]{$y$}}
\put(-0.1,-0.1){\makebox(0,0)[c]{0}}
\put(-0.1,1){\makebox(0,0)[c]{1}}
\put(1,-0.1){\makebox(0,0)[c]{1}}
\end{picture}
\end{center}
\caption{Geometry for the lid-driven cavity problem.}
\label{drcv}
\end{figure}

The problem consists of calculating the flow in a square cavity with
uniformly moving lid. The domain and
boundary conditions of the lid-driven cavity problem are shown in
fig.~\ref{drcv}, where $u$ and $v$ denote the velocity in $x$- and
$y$-direction, respectively.

The equations are given by
\begin{equation}
\label{eq_NS}
\left. \begin{array}{rcl} - \mathbf{u} \cdot \nabla \mathbf{u} + \frac 1 {Re} \Delta 
\mathbf{u} - \nabla p & = & 0 ~, \\
\nabla \cdot \mathbf{u} & = & 0 ~. \end{array} \right\}\end{equation} 
For the discretization we use a symmetry-preserving space discretization \cite{Verstappen},
 which is stable and does not introduce artificial diffusion.
 Furthermore, the grid is stretched towards the boundaries in order to resolve the boundary
 layers. The ratio between largest and smallest mesh size is about 5. This also means that
 we really need to change to fluxes through grid cell boundaries instead of velocities in order
 to get the required property that all elements in $B$ have the same magnitude (see the
 beginning of section~\ref{sec:ilu}).
The convergence tolerance is set to $10^{-6}$ in these experiments.
 The system matrix is the Jacobian from the first step of the Newton method at
 the current Reynolds number. In order to avoid convergence problems of Newton's method, we
 use the result at the previous Reynolds-number as a starting solution (The Reynolds
 numbers used are shown in table \ref{ta_NavStok2D_1}).

We first focus on the effect of increasing the Reynolds-number
(cf. table \ref{ta_NavStok2D_1}). The convergence is not independent of the
Reynolds-number. In our view this is not surprising, because the underlying continuous 
problem changes with the Reynolds number and more and more eigenvalues are getting close to 
the origin. This is different from the dependence on the mesh, where the continuous problem 
stays the same and all eigenvalues near the origin stay at their place. 

Next we refine the grid at a high Reynolds-number of 8000, close to the
point (cf. \citeasnoun{Tiesinga_2002}) where the steady state becomes unstable; results are 
shown in table
\ref{ta_NavStok2D_2}.
Note that the number of iterations is going down as we decrease the mesh-size. This is
because with decreasing mesh-size the physical size of the subdomains is
decreasing if we keep the number of unknowns per subdomain the same.
As the physical subdomain decreases, the diffusion plays a more
important role than the advection on that scale. Since the approximations take place at the
subdomain scale, the convergence behavior tends to that of the Stokes problem.

\begin{table}[!ht]
\begin{tabular}{|c|r|r|r|r|r|r|r|r|}
\hline
Re & N & nnz & $\mathrm{N_S}$ & n & iter & fill 1 & fill 2 \\
\hline
 500 &  785 408 & 6 794 252 &  129 025 & 40 069 & 59 & 6.41 & 2.59 \\
 1000 &  785 408 & 6 794 252 &  129 025 & 40 069 & 73 & 6.39 & 2.59 \\
 2000 &  785 408 & 6 794 252 &  129 025 & 40 069 & 87 & 6.38 & 2.65 \\
 4000 &  785 408 & 6 794 252 &  129 025 & 40 069 &  104 & 6.35 & 2.78 \\
 8000 &  785 408 & 6 794 252 &  129 025 & 40 069 &  130 & 6.33 & 2.72 \\
\hline
\end{tabular}
\caption{2D Driven cavity - increasing Reynolds-number, grid-size
$n_x=512$\label{ta_NavStok2D_1}}
\end{table}

\begin{table}[!ht]
\begin{tabular}{|c|r|r|r|r|r|r|r|r|}
\hline
$n_x$ & N & nnz & $\mathrm{N_S}$ & n & iter & fill 1 & fill 2 \\
\hline
64 & 12 160 &  103 820 & 1 793 &  533 &  185 & 6.09 & 0.418 \\
 128 & 48 896 &  420 620 & 7 681 & 2 341 &  181 & 6.22 & 0.953 \\
 256 &  196 096 & 1 693 196 & 31 745 & 9 797 &  167 & 6.29 & 1.75 \\
 512 &  785 408 & 6 794 252 &  129 025 & 40 069 &  130 & 6.33 & 2.72 \\
\hline
\end{tabular}
\caption{2D Driven cavity - grid refinement at $Re=8 000$\label{ta_NavStok2D_2}}
\end{table}

We conclude by mentioning that with the resulting preconditioner it was also quite easy to 
compute eigenvalues using MATLAB's eigs routine (i.e. ARPACK). Hence we can now study the 
stability problem near the point where the steady state becomes unstable using eigenvalue 
analysis.

\section{Discussion and conclusions}\label{sec:Concl}

In this paper we have shown that the structure preserving complete
$LDL^T$ factorization introduced in \citeasnoun{de_Niet_Wubs_2009}
of an \Fmat~can be transformed into an incomplete factorization.
We constructed an iterative solver for the whole system, which avoids having to balance
inner and outer iterations as in a segregated approach. Depending only on a single parameter
(the subdomain size), the method is as easy to use as a direct solver
and gives reliable results in a reasonable turn-around time. \par

For Stokes matrices we were able to prove grid-independent convergence. The total number of
operations required is currently not grid-independent since we use a direct solver to solve
the reduced system. However, the amount of work required for this step
is reduced by about the cube of the subdomain size in 2D and the sixth power
in 3D. So increasing the subdomain size by a factor 2 means in 2D a factor 8 and in 3D a factor 64.
For the Navier-Stokes equations we also observed grid-independent convergence.
We are developing a parallel C++ implementation of the method that can be applied
recursively, making it a multi-level method.\par
We proved the robustness of the method for Stokes and Navier-Stokes equations,
 where in the latter case the matrix should be definite. Computations show that the method
 still performs well for cases where eigenvalues pass the imaginary axis away from the
 origin (Hopf bifurcations).\par

In the case of \Fmats~ we are able to keep the computation in the kernel of the constraint
equation, i.e. for Stokes in the divergence free space, allowing us to use the CG method.
Though the \Fmats~ seem to be a limited class due to the constraints on the sparsity pattern
in $B$, many applications lead to matrices of this type.



\section*{Acknowledgements}

Part of this work was done during a sabbatical leave of the first author
to the Numerical Analysis Group at Rutherford Appleton Laboratory in Chilton (UK).
We kindly thank that group for the hospitality and good-fellowship.
The research of the second author was funded by the Netherlands Organization for Scientific
Research, NWO, through the contract ALW854.00.028.

\end{document}